\documentclass[11pt, a4paper]{article}

\usepackage{amsmath,amssymb,amsthm,mathrsfs}
\usepackage{geometry}
\usepackage[utf8]{inputenc}
\usepackage[T1]{fontenc}
\usepackage{hyperref}
\usepackage{cite}

\hypersetup{
    colorlinks=true,
    linkcolor=blue,
    citecolor=red,
    urlcolor=blue
}

\newtheorem{theorem}{Theorem}[section]

\theoremstyle{remark}
\newtheorem{remark}[theorem]{Remark}

\newcommand{\Z}{\mathbb{Z}}

\title{\textbf{Explicit Lower Bounds for Dirichlet Series of Higher Power Representation Functions}}
\author{Mahipal Gurram}
\date{\today}

\begin{document}
\maketitle

\begin{abstract}
We investigate Dirichlet-type series generated by representation functions that count the number of ways an integer can be expressed as a sum of $k$ signed higher even powers. By combining generalized theta generating functions with a family of generalized cotangent series introduced in previous work, we derive two distinct explicit lower bounds for these series. The first estimate arises from a geometric restriction of the lattice to its diagonal, while the second utilizes H\"older's inequality on the integral representation of the series. The methods presented here avoid modular techniques and offer a flexible analytic framework for higher-power representation problems.
\end{abstract}

\section{Introduction}

Representation problems in number theory, which ask in how many ways an integer can be expressed as a sum of powers, occupy a central position in both classical and modern analytic number theory. The case of sums of squares is particularly well understood and is deeply connected with modular forms, Jacobi theta functions, and special values of the Riemann zeta function \cite{HardyWright2008, Apostol1976}. Beyond quadratic forms, however, analogous questions involving higher even powers exhibit substantially more analytic flexibility and lack the rigid modular structure that governs the classical theory.

Let $m,k\in\mathbb{N}$. In this paper, we study the representation function
\[
r_{m,k}(n)
=
\#\Big\{(x_1,\dots,x_k)\in\mathbb{Z}^k :
x_1^{2m}+\cdots+x_k^{2m}=n
\Big\},
\]
which counts the number of representations of a nonnegative integer $n$ as a sum of $k$ even powers of order $2m$. We investigate the associated Dirichlet-type series
\[
S_{m,k}(a) = \sum_{n=0}^{\infty}\frac{r_{m,k}(n)}{n+a},
\qquad a>0.
\]
A necessary condition for the convergence of this series is $k < 2m$. Indeed, asymptotic estimates for the number of lattice points in $l^{2m}$-balls imply that $r_{m,k}(n) \sim C n^{\frac{k}{2m}-1}$ on average. Consequently, the general term of the series behaves like $n^{\frac{k}{2m}-2}$, and convergence requires $2 - \frac{k}{2m} > 1$, or $k < 2m$.

The analytic viewpoint adopted in this work is inspired by ideas traceable to Ramanujan and later Andrews \cite{Andrews1976}, utilizing generating functions as analytic objects capable of encoding deep arithmetic structure.The present paper builds on the author’s recent work
\cite{gurram2025generalizedcotangentserieslinks}, where the generalized cotangent
series
\[
U_{2n}(z)=\sum_{k\in\mathbb{Z}}\frac{1}{k^{2n}+z^{2n}}
\]
was introduced and shown to admit explicit finite trigonometric–hyperbolic
representations as well as an integral representation involving generalized
theta functions. These identities allow higher-dimensional lattice sums to be
reduced to one-dimensional analytic expressions.

A key observation is that the generalized theta function
\[
\Theta_m(q)=\sum_{k\in\mathbb{Z}} q^{k^{2m}}, \qquad 0<q<1,
\]
generates the representation numbers via
\[
\Theta_m(q)^k=\sum_{n=0}^{\infty} r_{m,k}(n)\,q^n.
\]
This identity enables the combination of geometric lattice restrictions with
analytic inequalities to obtain explicit bounds for the associated Dirichlet
series.In this paper, we utilize the relationship between $U_{2m}(z)$ and generalized theta functions to derive two explicit lower bounds for $S_{m,k}(a)$: one based on lattice geometry and one based on analytic convexity (H\"older's inequality).

\section{The Geometric Lower Bound}

Our first result relies on the geometric structure of the underlying lattice $\Z^k$.

\begin{theorem}[Geometric Lower Bound]\label{thm:geometric}
Let $m,k\in\mathbb{N}$ such that $k < 2m$, and let $a>0$. Let $U_{2m}(z)$ be the generalized cotangent series defined by
\[
U_{2m}(z)
=
\sum_{t\in\mathbb{Z}}\frac{1}{t^{2m}+z^{2m}},
\qquad z>0.
\]
Then the Dirichlet series associated with $r_{m,k}(n)$ satisfies
\[
\sum_{n=0}^{\infty}\frac{r_{m,k}(n)}{n+a}
\;\ge\;
\frac{1}{k}\,
U_{2m}\!\left(\left(\tfrac{a}{k}\right)^{1/(2m)}\right).
\]
\end{theorem}

\begin{proof}
By definition,
\[
\sum_{n=0}^{\infty}\frac{r_{m,k}(n)}{n+a}
=
\sum_{x_1,\dots,x_k\in\mathbb{Z}}
\frac{1}{x_1^{2m}+\cdots+x_k^{2m}+a}.
\]
Since the summand is strictly positive, we may obtain a lower bound by restricting the summation to any non-empty subset of the lattice $\mathbb{Z}^k$. We consider the diagonal sublattice defined by
\[
\mathcal{D} = \{(t,t,\dots,t)\in\mathbb{Z}^k : t\in\mathbb{Z}\}.
\]
Restricting the sum to $(x_1, \dots, x_k) \in \mathcal{D}$, we substitute $x_i = t$ for all $i=1,\dots,k$:
\[
\sum_{x_1,\dots,x_k\in\mathbb{Z}}
\frac{1}{x_1^{2m}+\cdots+x_k^{2m}+a}
\ge
\sum_{t\in\mathbb{Z}}
\frac{1}{k t^{2m}+a}.
\]
Factoring $k$ from the denominator yields
\[
\sum_{t\in\mathbb{Z}}
\frac{1}{k(t^{2m}+a/k)}
=
\frac{1}{k}
\sum_{t\in\mathbb{Z}}
\frac{1}{t^{2m}+(a/k)}.
\]
Recognizing the sum as $U_{2m}(z)$ with $z^{2m} = a/k$, we obtain
\[
\frac{1}{k}\,
U_{2m}\!\left(\left(\tfrac{a}{k}\right)^{1/(2m)}\right),
\]
which completes the proof.
\end{proof}

\section{The Analytic Lower Bound}

Our second result employs the integral representation of the Dirichlet series and H\"older's inequality.

\begin{theorem}[Analytic Lower Bound]\label{thm:holder}
Let $m,k\in\mathbb{N}$ with $1 \le k < 2m$, and let $a>0$. Then
\[
\sum_{n=0}^{\infty}\frac{r_{m,k}(n)}{n+a}
\;\ge\;
a^{\,k-1}
\left[
U_{2m}\!\left(a^{1/(2m)}\right)
\right]^k.
\]
\end{theorem}

\begin{proof}
Let $\Theta_m(q) = \sum_{t\in\mathbb{Z}} q^{t^{2m}}$ for $0<q<1$. Using the identity $\frac{1}{N+a} = \int_0^1 q^{N+a-1}dq$, we write the series as
\[
S_{m,k}(a) = \int_0^1 q^{a-1}\Theta_m(q)^k\,dq.
\]
We introduce the measure $d\mu(q) = q^{a-1}dq$ on $(0,1)$. The total measure is $\mu((0,1)) = \int_0^1 q^{a-1}dq = \frac{1}{a}$.
We apply H\"older's inequality in the form
\[
\int f \, d\mu \le \left(\int f^k \, d\mu \right)^{1/k} \left(\int 1^{k/(k-1)} \, d\mu \right)^{(k-1)/k},
\]
setting $f(q) = \Theta_m(q)$. This gives:
\[
\int_0^1 \Theta_m(q) q^{a-1}dq
\le
\left( \int_0^1 \Theta_m(q)^k q^{a-1}dq \right)^{1/k}
\left( \frac{1}{a} \right)^{(k-1)/k}.
\]
The term on the left is exactly $U_{2m}(a^{1/(2m)})$. The integral on the right is $S_{m,k}(a)$. Rearranging, we get:
\[
U_{2m}(a^{1/(2m)}) \le S_{m,k}(a)^{1/k} \cdot a^{-(k-1)/k}.
\]
Raising both sides to the power $k$ and solving for $S_{m,k}(a)$ yields:
\[
S_{m,k}(a) \ge a^{k-1} \left[ U_{2m}(a^{1/(2m)}) \right]^k.
\]
\end{proof}
\section{Illustrative Examples}

To demonstrate the explicitness of the bounds derived in Theorems \ref{thm:geometric},\ref{thm:holder} and we will use some result from this paper\cite{gurram2025generalizedcotangentserieslinks}, we examine the cases $m=1$ (quadratic forms) and $m=2$ (quartic forms).

\subsection{Case $m=1$: Sums of Squares}
In this case, $2m=2$. The generalized cotangent series reduces to the classical hyperbolic cotangent identity:
\[
U_2(z) = \frac{\pi}{z} \coth(\pi z).
\]
Although the strict convergence condition $k < 2m$ implies only $k=1$ is permissible here for the infinite Dirichlet series, the structural form of the bounds remains instructive.

\begin{itemize}
    \item \textbf{Geometric Bound:}
    \[
    \mathcal{B}_{\mathrm{geo}}(a) = \frac{1}{k} \left[ \frac{\pi}{\sqrt{a/k}} \coth\left(\pi \sqrt{\frac{a}{k}}\right) \right] = \frac{\pi}{\sqrt{ak}} \coth\left(\pi \sqrt{\frac{a}{k}}\right).
    \]
    
    \item \textbf{Analytic Bound:}
    \[
    \mathcal{B}_{\mathrm{ana}}(a) = a^{k-1} \left[ \frac{\pi}{\sqrt{a}} \coth(\pi \sqrt{a}) \right]^k = \frac{\pi^k}{\sqrt{a}^{2-k}} \coth^k(\pi \sqrt{a}).
    \]
\end{itemize}

\subsection{Case $m=2$: Sums of Fourth Powers}
Here $2m=4$. The condition $k < 2m$ allows for dimensions $k=1, 2, 3$. As shown in \cite{gurram2025generalizedcotangentserieslinks}, the series $U_4(z)$ admits the closed form:
\[
U_4(z) = \frac{\pi}{\sqrt{2}z^3} \left( \frac{\sin(\sqrt{2}\pi z) + \sinh(\sqrt{2}\pi z)}{\cosh(\sqrt{2}\pi z) - \cos(\sqrt{2}\pi z)} \right).
\]
Let $\Psi(x) = \frac{\sin(x) + \sinh(x)}{\cosh(x) - \cos(x)}$. The bounds become:

\begin{itemize}
    \item \textbf{Geometric Bound:}
    \[
    \mathcal{B}_{\mathrm{geo}}(a) = \frac{\pi k^{3/4}}{\sqrt{2} a^{3/4}} \Psi\left(\sqrt{2}\pi (a/k)^{1/4}\right).
    \]
    \item \textbf{Analytic Bound:}
    \[
    \mathcal{B}_{\mathrm{ana}}(a) = a^{k-1} \left[ \frac{\pi}{\sqrt{2} a^{3/4}} \Psi(\sqrt{2}\pi a^{1/4}) \right]^k = \frac{\pi^k}{2^{k/2}} a^{\frac{k}{4}-1} \left[ \Psi(\sqrt{2}\pi a^{1/4}) \right]^k.
    \]
\end{itemize}
These expressions provide fully explicit, non-trivial lower estimates for the spectral sums of the operator related to $x_1^4 + \dots + x_k^4$.

\section{Detailed Comparative Analysis}

We now compare the Geometric Lower Bound ($B_{\text{geo}}$) and the Analytic Lower Bound ($B_{\text{ana}}$). This comparison reveals a fundamental dichotomy between the geometric constraints of the lattice and the convexity properties of the generating functions.

\subsection{Dependence on the Spectral Parameter $a$}
As derived in the asymptotic analysis, the ratio of the bounds scales as:
\[
\mathcal{R}(a) := \frac{\mathcal{B}_{\mathrm{ana}}(a)}{\mathcal{B}_{\mathrm{geo}}(a)} \propto a^{\frac{k-1}{2m}}.
\]
Since $k \ge 1$, the exponent is non-negative.
\begin{enumerate}
    \item \textbf{The Singular Limit ($a \to 0^+$):} Since $\lim_{a \to 0} \mathcal{R}(a) = 0$, the \textbf{Geometric Bound} is significantly stronger. Physically, as the mass parameter $a$ vanishes, the sum is dominated by the lowest-lying energy states (lattice points near the origin). The diagonal restriction $\mathcal{D}$ captures these low-lying points efficiently ($1, 1, \dots, 1$), whereas the integral averaging in the H\"older bound smooths over them, resulting in a weaker estimate near the singularity.
    
    \item \textbf{The Large Shift Limit ($a \to \infty$):} Since $\lim_{a \to \infty} \mathcal{R}(a) = \infty$, the \textbf{Analytic Bound} dominates. For large $a$, the discrete granularity of the lattice becomes less significant than the global volume growth. The convexity inherent in H\"older's inequality captures this aggregate behavior more effectively than a restriction to a 1D subspace.
\end{enumerate}

\subsection{Dimensional Stability vs. Convexity}
A crucial distinction arises when considering the dependence on the dimension $k$ (within the convergent range $1 \le k < 2m$).

The two bounds arise from fundamentally different principles, leading to contrasting behaviors in $k$:

\begin{enumerate}
    \item \textbf{Polynomial Decay of the Geometric Bound:} 
    The geometric bound $\mathcal{B}_{\mathrm{geo}} \approx \frac{1}{k} U_{2m}((a/k)^{1/2m})$ scales with a polynomial prefactor $1/k$. This implies that the bound remains relatively \textbf{stable in high dimensions}. It does not degrade rapidly as $k$ increases; it merely scales down linearly to account for the density of the diagonal sublattice relative to the whole space.
    
    \item \textbf{Exponential Growth of the Analytic Bound:}
    The analytic bound $\mathcal{B}_{\mathrm{ana}} \approx [U_{2m}]^k$ exhibits \textbf{exponential dependence} on $k$. 
    \begin{itemize}
        \item When $U_{2m}$ is large (i.e., small $a$), this exponential growth makes the analytic bound extremely sensitive to dimension. However, due to the H\"older inequality direction, this "growth" actually represents a looser lower bound in the limit $a \to 0$ compared to the geometric case.
        \item Conversely, for fixed $a$ and $k$, this exponential non-linearity allows the bound to adapt to the "curvature" of the solution space better than the linear geometric restriction.
    \end{itemize}
\end{enumerate}

\begin{remark}[Regime Transition]
This analysis reveals a transition between a \textbf{geometry-dominated regime} and a \textbf{convexity-dominated regime}. 
\begin{itemize}
    \item When $a$ is small, the discrete lattice geometry dictates the behavior, and the polynomial stability of the diagonal restriction yields the optimal result.
    \item When $a$ is large, the analytic convexity dictates the behavior, and the exponential power of the H\"older inequality yields the optimal result.
\end{itemize}
\end{remark}

\section{Conclusion}
We have established two explicit lower bounds for Dirichlet series of higher power representation functions. The results highlight a necessary convergence condition $k < 2m$ and demonstrate that while analytic techniques (H\"older) provide elegant closed forms, simple geometric lattice restrictions can yield asymptotically stronger estimates near singularities.

\bibliographystyle{plain}
\bibliography{refs}
\end{document}